\documentclass[12pt,reqno]{amsart}
\usepackage{latexsym,amssymb,amscd}
\def\frk{\frak}               % font for "Fraktur"

\def\kk{{\Bbbk}}
\def\B'c{{\mathcal{B'}}}
\def\U'c{{\mathcal{U'}}}

\def\opn#1#2{\def#1{\operatorname{#2}}} % to make operators
\opn\chara{char}
\opn\length{\ell}
\opn\projdim{proj\,dim}
\opn\injdim{inj\,dim}
\opn\ini{in}
\opn\rank{rank}
\opn\Tiefe{Tiefe}
\opn\grade{grade}
\opn\height{height}
\opn\embdim{emb\,dim}
\opn\codim{codim}

\opn\Tr{Tr}
\opn\bigrank{big\,rank}
\opn\superheight{superheight}\opn\lcm{lcm}
\opn\trdeg{tr\,deg}%
\opn\reg{reg}
\opn\lreg{lreg}
\opn\deg{deg}
\opn\lcm{lcm}
\opn\syz{syz}

\opn\div{div}
\opn\Div{Div}
\opn\cl{cl}
\opn\Cl{Cl}
\opn\Spec{Spec}
\opn\Supp{Supp}
\opn\supp{supp}
\opn\Sing{Sing}
\opn\Ass{Ass}
\opn\Min{Min}
\opn\Ann{Ann}
\opn\Rad{Rad}
\opn\Soc{Soc}
\opn\Ker{Ker}
\opn\Coker{Coker}
\opn\Im{Im}
\opn\Hom{Hom}
\opn\Tor{Tor}
\opn\Ext{Ext}
\opn\End{End}
\opn\Aut{Aut}
\opn\id{id}

\opn\nat{nat}
\opn\GL{GL}
\opn\SL{SL}
\opn\mod{mod}
\opn\ord{ord}
\opn\depth{depth}
\opn\set{set}
\opn\Shad{Shad}
\opn\aff{aff}
\opn\con{conv}
\opn\relint{relint}
\opn\st{st}
\opn\lk{lk}
\opn\cn{cn}
\opn\core{core}
\opn\vol{vol}
\opn\Cut{Cut}
\opn\Mon{Mon}
\opn\gr{gr}

\def\pot#1#2{#1[\kern-0.28ex[#2]\kern-0.28ex]}

\opn\dirlim{\underrightarrow{\lim}}
\opn\invlim{\underleftarrow{\lim}}
\def\pnt{{\raise0.5mm\hbox{\large\bf.}}}

\def\twoline#1#2{\aoverb{\scriptstyle {#1}}{\scriptstyle {#2}}}
\newcommand{\aoverb}[2]{{\genfrac{}{}{0pt}{1}{#1}{#2}}}
\def\Implies{\ifmmode\Longrightarrow \else
     \unskip${}\Longrightarrow{}$\ignorespaces\fi}
\def\implies{\ifmmode\Rightarrow \else
     \unskip${}\Rightarrow{}$\ignorespaces\fi}
\def\iff{\ifmmode\Longleftrightarrow \else
     \unskip${}\Longleftrightarrow{}$\ignorespaces\fi}

\let\:=\colon
\newtheorem{Theorem}{Theorem}[section]
\newtheorem{Lemma}[Theorem]{Lemma}
\newtheorem{Corollary}[Theorem]{Corollary}
\newtheorem{Proposition}[Theorem]{Proposition}
\newtheorem{Remark}[Theorem]{Remark}

\let\epsilon=\varepsilon
\let\phi=\varphi
\let\kappa=\varkappa
\title{Monomial cut ideals}
\author{Anda Olteanu}
\address{Faculty of Mathematics and Computer Science, Ovidius University, Bd.\ Mamaia 124,
 900527 Constanta, Romania,} \email{olteanuandageorgiana@gmail.com}
 \thanks{The author was supported by the CNCSIS-UEFISCDI project PN II-RU PD 23/06.08.2010 and by the strategic grant POSDRU/89/1.5/S/58852, Project ``Postdoctoral program for training scientific researchers" co-financed by the European Social Fund within the Sectorial Operational Program Human Resources Development 2007 - 2013}

\begin{document}

\maketitle

\begin{abstract} B. Sturmfels and S. Sullivant associated to any graph a toric ideal, called the cut ideal. We consider monomial cut ideals and we show that their algebraic properties such as the minimal primary decomposition, the property of having a linear resolution or being Cohen--Macaulay may be derived from the combinatorial structure of the graph. 



\end{abstract}

\section*{Introduction} 
Defined by B. Sturmfels and S. Sullivant in \cite{SS}, cut ideals are generalizations of certain classes of toric ideals which appear in phylogenetics and in algebraic statistics. To any finite simple graph $G$ with the vertex set $V(G)$ and the set of edges $E(G)$, they associate a toric ideal, called \textit{the cut ideal of the graph $G$}, and denoted by $I_G$, as follows: Given a partition $A|B$ of the vertex set $V(G)$, one may consider the set $\Cut(A|B):=\{\{i,j\}\in E(G): i\in A,j\in B\mbox{ or }i\in B,j\in A\}$. Moreover, we assume that the partition $A|B$ is unordered. We consider the polynomial rings $\kk[q]:=\kk[q_{A|B}:A\cup B=V(G),\ A\cap B=\emptyset]$ and $S=\kk[s_{ij},t_{ij}:\{i,j\}\in E(G)]$ over a field $\kk$, and let $\Phi_G$ be the homomorphism of polynomial rings $\Phi_G:\kk[q]\rightarrow S$ defined by $$\Phi_G(q_{A|B})=u_{A|B}:=\prod_{\{i,j\}\in\Cut(A,B)}s_{ij}\prod_{\{i,j\}\in E(G)\setminus\Cut(A,B)}t_{ij}.$$ The \textit{cut ideal of $G$} is $I_G=\ker(\Phi_G)$.

Cut ideals have been intensively studied \cite{E}, \cite{NP}, \cite{SS}, but their algebraic properties are largely unknown. 
 
In this paper, we consider \textit{the monomial cut ideal}, that is the monomial ideal generated by the monomials which generate the toric ring $\Im(\Phi_G)$. More precisely, for a finite, simple, connected graph $G=(V(G),E(G))$, we consider the monomial ideal $I(G)=(u_{A|B}:A\cup B=V(G),A\cap B=\emptyset)$, where $u_{A|B}$ is the monomial defined above. We aim to see how the algebraic and homological properties of the monomial cut ideal of $G$ are related to the combinatorial properties of the graph.

The paper is structured as follows: The first section is devoted to the introduction of monomial cut ideals. We aim to determine the relation between the operation of deleting edges from a graph and the behaviour of monomial cut ideals. These results will be intensively used in the next sections. 

In the second section, we consider the monomial cut ideals of trees. We compute the minimal primary decomposition, which will allow us to derive some invariants such as the depth and the Castelnuovo--Mumford regularity. Our result shows that the minimal primary decomposition can be written just by looking at the graph. In the end of this section, we compute the Betti numbers of monomial cut ideals associated to trees. Our formula for the Betti numbers depends only on the number of edges of the tree. 

Next, in the third section, we pay attention to cycles. In order to determine the minimal primary decomposition, we distinguish between odd cycles and even cycles. We show that the minimal primary decomposition can be easily written just by knowing the length of the cycle. We compute the Castelnuovo--Mumford regularity and the depth, and we show that the monomial cut ideals of cycles have a pure resolution.

The last section is devoted to the study of monomial cut ideals associated to arbitrary graphs. For an arbitrary graph, the minimal primary decomposition is described in terms of the minimal prime ideals of the monomial cut ideals associated to its subgraphs which are cycles. Therefore, monomial cut ideals ``recognize" all the cycles from the associated graph. We also pay attention to several algebraic properties of monomial cut ideals such as having a linear resolution, being Gorenstein, Cohen--Macaulay or sequentially Cohen--Macaulay. We characterize all the monomial cut ideals which have one of these properties. Our characterization is expressed in terms of the combinatorial structure of the graph. 

\section{Monomial cut ideals}
Through this paper, by a graph $G=(V(G),E(G))$ we will understand a finite, simple, connected graph with the vertex set $V(G)$ and the set of edges $E(G)$. We denote by $\mathcal{P}(V(G))$ the set of all the unordered partitions of $V(G)$. For $A|B\in\mathcal{P}(V(G))$, one may consider the set $$\Cut(A|B):=\{\{i,j\}\in E(G): i\in A,j\in B\mbox{ or }i\in B,j\in A\}.$$ We consider the polynomial ring $S=\kk[s_{ij},t_{ij}:\{i,j\}\in E(G)]$ over a field $\kk$. To every $A|B\in\mathcal{P}(V(G))$, one may associate a squarefree monomial $$u_{A|B}:=\prod_{\{i,j\}\in\Cut(A,B)}s_{ij}\prod_{\{i,j\}\in E(G)\setminus\Cut(A,B)}t_{ij}.$$ The ideal $I(G)=(u_{A|B}:A|B\in\mathcal{P}(V(G)))$ will be called \textit{the monomial cut ideal of the graph $G$}. 

For a monomial $m$ in $S$, we denote $\supp(m)=\{s_{ij}: s_{ij}\mid m\}\cup\{t_{ij}: t_{ij}\mid m\}$. Note that, any monomial $m$ in $S$ can be uniquely written  as $m=m_sm_t$, where $\supp(m_s)=\supp(m)\cap\{s_{ij}:\{i,j\}\in E(G)\}$ and  $\supp(m_t)=\supp(m)\cap\{t_{ij}:\{i,j\}\in E(G)\}$. We consider $\deg_{s}(m)=\deg(m_s)$ the $s$-degree of $m$, $\deg_{t}(m)=\deg(m_t)$ the $t$-degree of $m$, and $\deg(m)=\deg_s(m)+\deg_t(m)$ the total degree of $m$. 
For an integer $d\geq2$, we will denote by $\Mon_d(S)$ the set of all squarefree monomials of degree $d$ in $S$. Given a monomial ideal $I$ in $S$, we denote by $\mathcal{G}(I)$ the minimal monomial set of generators.

One may note that the monomial cut ideal of the graph $G=(V(G),E(G))$ is generated by $2^{|V(G)|-1}$ squarefree monomials of degree $|E(G)|$, since the partitions are unordered. Moreover, for every $A|B\in\mathcal{P}(V(G))$ and every $\{i,j\}\in E(G)$, $\deg(\gcd(u_{A|B},s_{ij}t_{ij}))=1$.

\begin{Remark}\label{subgraph}\rm If $G$ is a graph and $H\subseteq G$ is a subgraph, then $I(G)\subseteq I(H)$, since any partition of $V(G)$ induces a partition of $V(H)$.
\end{Remark}
We are interested to see what happens with the monomial cut ideal if one deletes an edge from the graph. First of all, we assume that the graph contains whiskers and we delete one whisker from the graph. We will consider that, when we delete a whisker, we also remove the free vertex from the vertex set.

\begin{Proposition}\label{whisker} Let $G=(V(G),E(G))$ be a graph and $\{\alpha,\beta\}\in E(G)$ be a whisker. If $H$ is the subgraph of $G$ obtained from $G$ by deleting the edge $\{\alpha,\beta\}$, then $$I(G)=I(H)\cap(s_{\alpha\beta},t_{\alpha\beta}).$$
\end{Proposition}

\begin{proof}``$\subseteq$'' Since $H\subseteq G$ is a subgraph, by Remark~\ref{subgraph}, we have that $I(G)\subset I(H)$. Also, $I(G)\subset(s_{\alpha\beta},t_{\alpha\beta})$, since $\{\alpha,\beta\}\in E(G)$.

``$\supseteq$'' Since $\gcd(s_{\alpha\beta}t_{\alpha\beta},m)=1$, for all $m\in \mathcal{G}(I(H))$, we have to prove that $ms_{\alpha\beta},\ mt_{\alpha\beta}\in I(G)$, for any monomial $m\in \mathcal{G}(I(H))$.

Since $\{\alpha,\beta\}\in E(G)$ is a whisker, we may assume that $\beta$ is the free vertex, hence $\beta\notin V(H)$. Therefore $V(H)=V(G)\setminus\{\beta\}$. Let $A|B\in\mathcal{P}(V(H))$ and $u_{A|B}\in \mathcal{G}(I(H))$. If $\alpha\in A$, then $u_{A|B}s_{\alpha\beta}\in I(G)$, since $A|(B\cup\{\beta\})\in\mathcal{P}(V(G))$. Also, $(A\cup\{\beta\})|B\in\mathcal{P}(V(G))$, therefore $u_{A|B}t_{\alpha\beta}=u_{A\cup\{\beta\}|B}\in I(G)$. If $\alpha\notin A$, then $u_{A|B}t_{\alpha\beta}\in I(G)$, since $A|(B\cup\{\beta\})\in\mathcal{P}(V(G))$. Also, $(A\cup\{\beta\})|B\in\mathcal{P}(V(G))$, therefore $u_{A|B}s_{\alpha\beta}=u_{A\cup\{\beta\}|B}\in I(G)$.
\end{proof}

\begin{Proposition}\label{edge} Let $G=(V(G),E(G))$ be a graph, $\{\alpha,\beta\}\in E(G)$, and $H$ be the subgraph of $G$ obtained from $G$ by deleting the edge $\{\alpha,\beta\}$. Then $$I(H)=I(G):(s_{\alpha\beta}t_{\alpha\beta}).$$
\end{Proposition}

\begin{proof} ``$\subseteq$'' Let $m\in \mathcal{G}(I(H))$, that is there exists $A|B\in\mathcal{P}(V(H))$ such that $$m=u_{A|B}=\prod_{\{i,j\}\in\Cut_H(A,B)}s_{ij}\prod_{\{i,j\}\in E(H)\setminus\Cut_H(A,B)}t_{ij}.$$ We have $\deg(m)=|E(G)|-1$.

If $\{\alpha,\beta\}$ is not a whisker, then $V(H)=V(G)$ and $A|B\in\mathcal{P}(V(G))$. Hence, we may write
$$m=u_{A|B}=\prod_{\twoline{\{i,j\}\in\Cut_G(A,B)}{\{i,j\}\neq\{\alpha,\beta\}}}s_{ij}\prod_{\twoline{\{i,j\}\in E(G)\setminus\Cut_G(A,B)}{\{i,j\}\neq\{\alpha,\beta\}}}t_{ij},$$ and $s_{\alpha\beta}t_{\alpha\beta}m\in I(G)$.

If $\{\alpha,\beta\}$ is a whisker, we may assume that $\beta$ is the free vertex, hence $V(G)=V(H)\cup\{\beta\}$, and $\beta\notin V(H)$. If $\alpha\in A$, then $A|(B\cup\{\beta\})\in\mathcal{P}(V(G))$. In this case $$s_{\alpha\beta}t_{\alpha\beta}m=t_{\alpha\beta}\prod_{\{i,j\}\in\Cut_G(A,B\cup\{\beta\})}s_{ij}\prod_{\{i,j\}\in E(G)\setminus\Cut_G(A,B\cup\{\beta\})}t_{ij}\in I(G).$$ If $\alpha\notin A$ $$s_{\alpha\beta}t_{\alpha\beta}m=s_{\alpha\beta}\prod_{\{i,j\}\in\Cut_G(A,B\cup\{\beta\})}s_{ij}\prod_{\{i,j\}\in E(G)\setminus\Cut_G(A,B\cup\{\beta\})}t_{ij}\in I(G).$$

``$\supseteq$'' Let $m\in I(G):(s_{\alpha\beta}t_{\alpha\beta})$ be a minimal monomial generator. Therefore, there exists $A|B\in\mathcal{P}(V(G))$ such that $ m=u_{A|B}/\gcd(u_{A|B},s_{\alpha\beta}t_{\alpha\beta})$, $\ u_{A|B}\in \mathcal{G}(I(G)).$
Obviously, $\deg(\gcd(u_{A|B},s_{\alpha\beta}t_{\alpha\beta}))=1$.

Firstly, we consider that $\gcd(u_{A|B},s_{\alpha\beta}t_{\alpha\beta})=s_{\alpha\beta}$, therefore $\{\alpha,\beta\}\in\Cut(A,B)$. We have

$$m=\frac{1}{s_{\alpha\beta}}\prod_{\{i,j\}\in\Cut_G(A,B)}s_{ij}\prod_{\{i,j\}\in E(G)\setminus\Cut_G(A,B)}t_{ij}=$$ $$=\prod_{\{i,j\}\in\Cut_G(A,B)\setminus\{\alpha,\beta\}}s_{ij}\prod_{\{i,j\}\in E(G)\setminus\Cut_G(A,B)}t_{ij}.$$

If $\{\alpha,\beta\}$ is not a whisker of $G$, then $A|B\in\mathcal{P}(V(H))$. In this case we get that $\Cut_H(A|B)=\Cut_G(A,B)\setminus\{\alpha,\beta\}$ and $E(H)\setminus\Cut_H(A,B)=E(G)\setminus \Cut_G(A,B),$ therefore $m\in I(H)$.

If $\{\alpha,\beta\}$ is a whisker of $G$, then we may assume that $\beta$ is the free vertex, and $\alpha\in A$ and $\beta\in B$. One may note that $A|(B\setminus\{\beta\})\in\mathcal{P}(V(H))$,
$\Cut_H(A,B\setminus\{\beta\})=\Cut_G(A,B)\setminus\{\alpha,\beta\}$ and $E(H)\setminus\Cut_H(A,B\setminus\{\beta\})=E(G)\setminus \Cut_G(A,B),$ therefore $m\in I(H)$.

Let us assume that $\gcd(u_{A|B},s_{\alpha\beta}t_{\alpha\beta})=t_{\alpha\beta}$, therefore $\{\alpha,\beta\}\notin\Cut(A,B)$. We have

$$m=\frac{1}{t_{\alpha\beta}}\prod_{\{i,j\}\in\Cut_G(A,B)}s_{ij}\prod_{\{i,j\}\in E(G)\setminus\Cut_G(A,B)}t_{ij}=$$ $$=\prod_{\{i,j\}\in\Cut_G(A,B)}s_{ij}\prod_{\{i,j\}\in (E(G)\setminus\{\alpha,\beta\})\setminus\Cut_G(A,B)}t_{ij}.$$

If $\{\alpha,\beta\}$ is not a whisker of $G$, then $A|B\in\mathcal{P}(V(H))$. In this case $\Cut_H(A|B)=\Cut_G(A,B)$, and $E(H)\setminus\Cut_H(A,B)=(E(G)\setminus\{\alpha,\beta\})\setminus \Cut_G(A,B),$ therefore $m\in I(H)$.

If $\{\alpha,\beta\}$ is a whisker of $G$, then we may assume that $\beta$ is the free vertex, and $\alpha,\beta\in A$. One may note that $(A\setminus\{\beta\})|B\in\mathcal{P}(V(H))$, $\Cut_H(A\setminus\{\beta\},B)=\Cut_G(A,B)$, and $E(H)\setminus\Cut_H(A\setminus\{\beta\},B)=(E(G)\setminus\{\alpha,\beta\})\setminus \Cut_G(A,B),$ therefore $m\in I(H)$.
\end{proof}
\section{Monomial cut ideals of trees}

In this section, we aim at determining the minimal primary decomposition of monomial cut ideals associated to trees and at computing some invariants. 

\begin{Proposition}\label{primdectree} Let $G=(V(G),E(G))$ be a tree. Then the minimal primary decomposition of $I(G)$ is $$I(G)=\bigcap _{\{i,j\}\in E(G)}(s_{ij},t_{ij}).$$
\end{Proposition}

\begin{proof} ``$\supseteq$'' It is easy to see that $I(G)\subset(s_{ij},t_{ij})$, for all $\{i,j\}\in E(G)$, and that they are minimal prime ideals of $I(G)$.

``$\subseteq$'' We use induction on the number of edges. If $|E(G)|=1$, then $I(G)=(s_{12},t_{12})$. We assume that the statement is true for any tree with at most $r$ edges, $r\geq1$. Let $G$ be a tree with $r+1$ edges and $\{\alpha,\beta\}$ be a whisker. Let $H$ be the subgraph of $G$ obtained by deleting the edge $\{\alpha,\beta\}$. Then $H$ is a tree with $r$ edges. By Proposition \ref{whisker} and by the induction hypothesis $$I(G)=I(H)\cap(s_{\alpha\beta},t_{\alpha\beta})=\left(\bigcap_{\{i,j\}\in E(H)}(s_{ij},t_{ij})\right)\cap(s_{\alpha\beta},t_{\alpha\beta})=\bigcap_{\{i,j\}\in E(G)}(s_{ij},t_{ij}),$$which ends the proof.
\end{proof}

\begin{Remark}\label{line}\rm By the minimal primary decomposition, one may easy note that, if $G=(V(G),E(G))$ and $G'=(V(G),E(G'))$ are two trees with the same set of vertices and $|E(G)|=|E(G')|$, then $S/I(G)$ and $S'/I(G')$ are isomorphic as $\kk$-algebras, where $S$ and $S'$ are the corresponding polynomial rings. Therefore, in order to compute algebraic and homological invariants of the monomial cut ideals of trees, it is enough to consider only the case of the path graph. 
\end{Remark}

In the rest of this section, we will denote by $T_r$ the path graph, on the set of vertices $\{1,\ldots,r+1\}$ and the set of edges $E(T_r)=\{\{i,i+1\}:1\leq i\leq r\}$. 

We consider the set of squarefree monomials $$M_r(T_r)=\{m\in\Mon_r(S):s_{i,i+1}t_{i,i+1}\nmid m,\mbox{ for all }1\leq i\leq r\}.$$

\begin{Remark}\rm If $T_r$ is the path graph, $r\geq1$, then $M_r(T_r)$ is the minimal monomial set of generators of $I(T_r)$.
\end{Remark}

 Firstly, we will show that the monomial cut ideal of the path graph has linear quotients. We recall that, if $I\subseteq\kk[x_1,\ldots,x_n]$ is a squarefree monomial ideal with $\mathcal{G}(I)=\{u_1,\ldots,$ $u_\ell\}$, then $I$ has linear quotients, \cite{H}, if and only if, all $i$ and for all $j<i$ there exist an integer $k<i$ and an integer $l$ such that
$$\frac{u_k}{\gcd(u_k,u_i)}=x_l\ \mbox{and}\ x_l\ \mbox{divides}\ \frac{u_j}{\gcd(u_j,u_i)}.$$ 
 \begin{Proposition} Let $T_r$ be the path graph, $r\geq1$. Then $I(T_r)$ has linear quotients.
 \end{Proposition} 
 \begin{proof} Let us fix on $S$ the lexicographical order with $s_{1,2}>s_{2,3}>\cdots>s_{r,r+1}>t_{1,2}>\cdots>t_{r,r+1}$ and let us denote $R=2^r$. We consider $\mathcal{G}(I(T_r))=\{u_1,\ldots,u_R\}$ with $u_1>_{lex}\cdots>_{lex}u_{R}$. Let $1\leq\alpha<\beta\leq R$. We have that $u_{\alpha}>_{lex}u_{\beta}$, hence there exists $s_{i,i+1}\in\supp(u_{\alpha})\setminus \supp(u_{\beta})$. One may note that the case when $s_{i,i+1}\notin\supp(u_{\alpha})\setminus\supp(u_{\beta})$, for all $\{i,i+1\}\in E(T_r)$, is impossible. Since $s_{i,i+1}\nmid u_{\beta}$, we must have $t_{i,i+1}\mid u_{\beta}$. Let us consider the monomial $u_{\gamma}=s_{i,i+1}u_{\beta}/t_{i,i+1}$. We obviously have $u_{\gamma}>_{lex}u_{\beta}$, thus $\gamma<\beta$, and $u_{\gamma}\in M_r(T_r)$ since $s_{j,j+1}t_{j,j+1}\nmid u_{\gamma}$ for any $1\leq j\leq r$ and $\deg(u_{\gamma})=r$. Moreover, $u_{\gamma}/\gcd(u_{\gamma},u_{\beta})=s_{i,i+1}$. 
 \end{proof}
 The above result has the following immediate consequence, \cite{HeTa}:
 \begin{Corollary}\label{linrestrees} Let $T_r$ be the path graph, $r\geq1$. Then $I(T_r)$ has a linear resolution.
 \end{Corollary}
 
 For a monomial ideal with linear quotients generated in one degree, we can compute the Betti numbers \cite{HeTa}, \cite[p. 135]{HeHi}. Firstly, we fix some more notations. Let $I$ be a monomial ideal of $S$ with $\mathcal{G}(I)=\{u_1,\ldots,u_\ell\}$ and assume that $I$ has linear quotients with respect to the sequence $u_1,\ldots,u_\ell$. We denote $\set(u_k)=\mathcal{G}((u_1,\ldots,u_{k-1}):u_k)$ and $r_k=|\set(u_k)|$, for all $2\leq k\leq \ell$. Then one has
\[\beta_i(I)=\sum_{k=2}^{\ell}\left(\twoline{r_k}{i} \right),
\]
for all $i$. Moreover, $\projdim(I)=\max\{r_k:2\leq k\leq \ell\}$. 

For monomial cut ideals of the path graph, the numbers $r_k$ can be determined. In the rest of this section, for a path graph $T_r$, we fix on $S$ the lexicographical order with $s_{1,2}>s_{2,3}>\cdots>s_{r,r+1}>t_{1,2}>\cdots>t_{r,r+1}$ and we denote $R=2^r$. We consider $\mathcal{G}(I(T_r))=\{u_1,\ldots,u_R\}$ with $u_1>_{lex}\cdots>_{lex}u_{R}$.

\begin{Proposition}\label{linquot} Let $T_r$ be the path graph, $r\geq1$. Then, in the above notations, $\set(u_k)=\{s_{i,i+1}:t_{i,i+1}\mid u_k\}$, $k\geq2$. In particular, $r_k=\deg_t(u_k)$.
\end{Proposition}
\begin{proof} ``$\supseteq$" One may note that, for every $k\geq2$, there exists some $i$, $1\leq i\leq r$, such that $t_{i,i+1}|u_{k}$. Let $k\geq2$ and $\{i,i+1\}\in E(T_r)$ such that $t_{i,i+1}\mid u_k$. Then the monomial $u_{\alpha}=s_{i,i+1}u_k/t_{i,i+1}\in M_r(T_r)=\mathcal{G}(I(T_r))$ and, since $u_{\alpha}>_{lex}u_k$, we have $\alpha<k$. Taking into account that $s_{i,i+1}=u_{\alpha}/\gcd(u_{\alpha},u_k)$, we have that $s_{i,i+1}\in\set(u_k)$.

``$\subseteq$" Firstly, we note that $t_{i,i+1}\notin \set(u_k)$, for every $\{i,i+1\}\in E(T_r)$. Indeed, if we assume that there exists $\{i,i+1\}\in E(T_r)$ such that $t_{i,i+1}\in\set(u_k)$, then there exists $\alpha<k$, such that $t_{i,i+1}=u_{\alpha}/\gcd(u_{\alpha},u_k)$. In particular, $\supp(u_{\alpha})=(\supp(u_k)\setminus\{s_{i,i+1}\})\cup\{t_{i,i+1}\}$, thus $u_{k}>_{lex}u_{\alpha}$ and $k<\alpha$, a contradiction. 

Let $s_{i,i+1}\in\set(u_k)$, that is there exists $\alpha<k$ such that $s_{i,i+1}=u_{\alpha}/\gcd(u_{\alpha},u_k)$. Thus, $s_{i,i+1}\notin\supp(u_k)$, therefore $t_{i,i+1}\mid u_k$, which ends the proof.
\end{proof}
We can determine the Betti numbers of monomial cut ideals of the path graph.

\begin{Proposition} Let $r\geq1$ be an integer and $T_r$ be the path graph. Then $$\beta_{i}(I(T_r))=2^{r-i}{{r}\choose{i}},$$for all $0\leq i\leq r$. Moreover, $\projdim(I(C_r))=r=E(C_r)$.
\end{Proposition}
 \begin{proof} Since $I(T_r)$ has linear quotients with respect to the decreasing lexicographical order of the minimal monomial generators, by Proposition~\ref{linquot} and taking into account that there are ${r}\choose{r_k}$ monomials of $t$-degree $r_k$ in $M_r(T_r)$, we have that
	\[\beta_{i}(I(T_r))=\sum_{k=1}^{r}{r\choose k}{k\choose{i}}=\sum_{k=1}^{r}{r\choose i}{r-i\choose{k-i}}={r\choose i}\sum_{k=1}^{r}{r-i\choose{k-i}}=2^{r-i}{{r}\choose{i}}.
\]In particular, $\projdim(I(C_r))=r$.
 \end{proof}
 Taking into account Remark~\ref{line}, we get the following result:
\begin{Corollary}\label{Bettitree} Let $G=(V(G),E(G))$ be a tree. Then $$\beta_{i}(I(G))=2^{|E(G)|-i}{{|E(G)|}\choose{i}},$$ for all $0\leq i\leq |E(G)|$.
\end{Corollary} 
In particular, we get the Catelnuovo--Mumford regularity and the depth of monomial cut ideals associated to trees.
\begin{Corollary} Let $G=(V(G),E(G))$ be a tree. Then 
\begin{itemize}
	\item[(i)] $\reg(S/I(G))=|E(G)|-1$.
	\item[(ii)] $\depth(S/I(G))=|E(G)|-1$.
\end{itemize}
\end{Corollary}
\begin{proof} 
(i) Using Remark~\ref{line} and Corollary~\ref{linrestrees}, one has that $I(G)$ is an ideal generated in degree $|E(G)|$ with a linear resolution, therefore the statement follows.

(ii)  The statement follows by Corollary~\ref{Bettitree}.
\end{proof}

\section{Monomial cut ideals of cycles}
Through this section, for a cycle $C_r=(V(C_r),E(C_r))$, we will consider that $V(C_r)=\{1,\ldots,r\}$ and we will denote the edges of $C_r$ by $\{i,i+1\}$ for $i=1,\ldots,r$. By convention, we will identify the ``edge" $\{r,r+1\}$ with the edge $\{1,r\}$. %We will also consider that the set of variables is $s_1,\ldots, s_r,t_1,\ldots, t_r$, where $s_i$ is the variable $s_{i,i+1}$ and $t_i$ is the variable $t_{i,i+1}$.

\begin{Proposition}\label{evensdeg} Let $C_r$ be a cycle, $r\geq3$ and $m\in I(C_r)$ be a minimal monomial generator. Then $\deg_s(m)$ is an even number.
\end{Proposition}
\begin{proof} 
One may note that every free vertex gives us a factor of $s$-degree $2$ and every sequence of adjacent vertices gives us a factor of $s$-degree $2$. Thus, $\deg_s(m)$ is an even number. Indeed, since $m\in \mathcal{G}(I(C_r))$, there exists $A|B\in\mathcal{P}(V(C_r))$ such that $$m=u_{A|B}=\prod_{\{i,i+1\}\in\Cut(A,B)}s_{i,i+1}\prod_{\{i,i+1\}\in E(G)\setminus\Cut(A,B)}t_{i,i+1}.$$ Let $i\in A$. If $i-1,i+1\notin A$, then $s_{i-1,i}s_{i,i+1}\mid m$. If $i,i+1,\ldots,i+j\in A$, for some $j\geq1$ and $i-1,i+j+1\notin A$, then $s_{i-1,i}s_{i+j,i+j+1}\mid m$.
\end{proof}
Let us fix some more notations. For a cycle $C_r$, $r\geq3$, we consider the following sets: 
$$M_e(C_r)=\{m\in \Mon_{r}(S)\ :\ \deg_s(m)=\mbox{even},\  s_{i,i+1}t_{i,i+1}\nmid m, \mbox{ for all}\ 1\leq i\leq r\}$$and $$M_o(C_r)=\{m\in \Mon_{r}(S)\ :\ \deg_s(m)=\mbox{odd},\  s_{i,i+1}t_{i,i+1}\nmid m, \mbox{ for all}\ 1\leq i\leq r\}.$$

For cycles, one may determine the minimal set of monomial generators of the corresponding monomial cut ideal.

\begin{Corollary}\label{gencycle} Let $r\geq3$ be an integer and $C_r$ a cycle. Then $M_e(C_r)$ is the minimal monomial set of generators of the monomial cut ideal $I(C_r)$.
\end{Corollary}

\begin{proof} By Proposition \ref{evensdeg}, $\mathcal{G}(I(C_r))\subseteq M_e(C_r)$. Since $|\mathcal{G}(I(C_r))|=2^{r-1}=|M_e(C_r)|$, the statement follows.
\end{proof}

Next, we aim at determining the minimal primary decomposition of monomial cut ideals of cycles. In order to do this, we have to distinguish between odd cycles and even cycles. We start with a remark concerning the height of the minimal prime ideals of the monomial cut ideals of cycles.

\begin{Lemma}\label{height2and|E|}Let $r\geq3$ be an integer and $C_r$ a cycle. Then
\begin{itemize}
	\item[(a)]  $\{\frk p\in\Min(S/I(C_r)):\height(\frk p)=2\}=\{(s_{i,i+1},t_{i,i+1}):1\leq i\leq r\}$.
	\item[(b)] If $\frk{p}\in\Min(S/I(C_r))$, $\height(\frk{p})\geq3$, then $\height(\frk{p})=|E(C_r)|=r$.
\end{itemize}
\end{Lemma}
\begin{proof} (a) Obviously, $(s_{i,i+1},t_{i,i+1})\supset I(C_r)$ for any $\{i,i+1\}\in E(C_r)$ and they are minimal prime ideals of $I(C_r)$. 

Let $\frk p\in\Min(S/I(C_r))$, $\height(\frk p)=2$. Let $\{i,i+1\},\{j,j+1\}\in E(C_r)$. Since $t_{1,2}\cdots t_{r,r+1}\in \mathcal{G}(I(C_r))$, we must have $t_{i,i+1}\in\frk p$ or $t_{j,j+1}\in\frk p$. One may note that $\frk{p}\neq(t_{i,i+1},t_{j,j+1})$. Indeed, if this is case, we can consider the monomial $m=s_{i,i+1}s_{j,j+1}\prod _{\{k,k+1\}\in E(C_r),\ k\neq i,j}t_{k,k+1}\in I(C_r)$, and $m\notin\frk p$, a contradiction. Let us assume that $t_{i,i+1}\in\frk p$. In this case, we may consider the monomial $m=s_{i,i+1}s_{k,k+1}\prod _{\{l,l+1\}\in E(C_r),\ l\neq i,k}t_{l,l+1}$, for some $k\neq j$. Since $m$ is divisible by $s_{i,i+1}t_{j,j+1}$ and $t_{j,j+1}\notin\frk p$, we must have $s_{i,i+1}\in\frk p$.

(b) If we assume that $\height(\frk p)>|E(C_r)|$, then there exists $\{i,i+1\}\in E(C_r)$ such that $(s_{i,i+1},t_{i,i+1})\subset \frk p$, contrary to $\frk p\in\Min(S/I(C_r))$. If we consider that $\height(\frk p)<|E(C_r)|$, then there exists a monomial $M\in \mathcal{G}(I(C_r))$ divisible by $\prod _{s_{i,i+1}\in\frk p}t_{i,i+1}\prod _{t_{i,i+1}\in\frk p}s_{i,i+1}$ and $M\notin\frk p$. Therefore, we must have $\height(\frk p)=|E(C_r)|$.
\end{proof}

We determine the minimal primary decomposition of the monomial cut ideals of odd cycles.

\begin{Proposition}\label{primdecodd} Let $r\geq3$ be an odd integer and $C_r$ be a cycle. Then the minimal primary decomposition of $I(C_r)$ is
	\[I(C_r)=\bigcap_{m\in M_e(C_r)}(\supp(m))\cap\bigcap_{\{i,i+1\}\in E(C_r)}(s_{i,i+1},t_{i,i+1}).
\]
\end{Proposition}

\begin{proof}``$\supseteq$" By Lemma~\ref{height2and|E|}, $(s_{i,i+1},t_{i,i+1})\supset I(C_r)$ for all $1\leq i\leq r$.

Let $m\in M_e(C_r)$ and $\frk{p}=(\supp(m))$. We have to show that $\frk{p}\in\Min(S/I(C_r))$. First of all, we will prove that $\frk{p}\supseteq I(C_r)$. Assume by contradiction that $\frk{p}\nsupseteq I(C_r)$, that is there exists a monomial $w\in \mathcal{G}(I(C_r))$ such that $w\notin\frk{p}$. We must have that $w\in M_e(C_r)$ since, $w\in \mathcal{G}(I(C_r))= M_e(C_r)$. In particular, $\deg_s(w)$ is an even number. 

On the other hand, $|\supp(m)|=|\supp(w)|=r$ and $\supp(w)\cap\supp(m)=\emptyset$. Thus $\supp(m)\cup\supp(w)=\{s_{i,i+1},t_{i,i+1}\ :\ \{i,i+1\}\in E(C_r)\}$. In particular, $w=\prod _{s_{i,i+1}\in\frk{p}}t_{i,i+1}\prod _{t_{i,i+1}\in\frk{p}}s_{i,i+1}$. Since $m\in M_e(C_r)$, we have $\deg_s(m)$ is an even number, thus $\deg_t(m)$ is an odd number, which implies $\deg_s(w)=\deg_t(m)$ is odd, a contradiction.

We prove that $\frk{p}\in\Min(S/I(C_r))$. Let $\{\alpha,\alpha+1\}\in E(C_r)$. If $s_{\alpha,\alpha+1}\in\frk p$, let $\frk{p}'=(\supp(m/s_{\alpha,\alpha+1}))$. Assume by contradiction that $\frk{p}'\supseteq I(C_r)$. We consider the monomial $M=s_{\alpha,\alpha+1}\prod _{s_{i,i+1}\in\frk{p}'}t_{i,i+1}\prod _{{t_{i,i+1}\in\frk{p}'}}s_{i,i+1}.$ Since $m\in M_e(C_r)$, $\deg_s(m)$ is an even number, which implies that $\deg_t(m)$ is an odd number. Since $\deg_s(M)=\deg_t(m)+1$, we get that $\deg_s(M)$ is an even number and $M\in \mathcal{G}(I(C_r))$, but $M\notin\frk{p}'$, a contradiction.

If $t_{\alpha,\alpha+1}\in\frk p$, let $\frk{p}'=(\supp(m/t_{\alpha,\alpha+1}))$. Assume by contradiction that $\frk{p}'\supseteq I(C_r)$. Let us consider the monomial $M=t_{\alpha,\alpha+1}\prod _{s_{i,i+1}\in\frk{p}'}t_{i,i+1}\prod _{{t_{i,i+1}\in\frk{p}'}}s_{i,i+1}$ which has $\deg_s(M)$ even, since $\deg_s(M)=\deg_t(m)-1$ and, taking into account that $\deg_s(m)$ is even, we have that $\deg_t(m)$ is an odd number. Therefore $M\in \mathcal{G}(I(C_r))$, but $M\notin\frk{p}'$, a contradiction.

``$\subseteq$''  By Lemma~\ref{height2and|E|}, $I(C_r)\subseteq\bigcap _{\{i,i+1\}\in E(C_r)}(s_{i,i+1},t_{i,i+1})$. Let $\frk p\in\Min(S/I(C_r))$, $\height(\frk p)\geq3$. By using Lemma~\ref{height2and|E|}, we must have $\height(\frk p)=|E(C_r)|$. Let us consider the monomial $m=\prod _{s_{i,i+1}\in\frk p}s_{i,i+1}\prod _{t_{i,i+1}\in\frk p}t_{i,i+1}.$ Assume by contradiction that $\deg_s(m)$ is odd. Since $C_r$ is an odd cycle, we must have $\deg_t(m)$ even. Therefore, the monomial $M=\prod _{s_{i,i+1}\in\frk p}t_{i,i+1}\prod _{t_{i,i+1}\in\frk p}s_{i,i+1}$ has $\deg_s(M)$ even, that is $M\in \mathcal{G}(I(C_r))$, but $M\notin\frk p$, a contradiction.
\end{proof}

The following result describes the minimal primary decomposition for monomial cut ideals of even cycles.
 
\begin{Proposition}\label{primdeceven} Let $r\geq3$ be an even integer and $C_r$ be a cycle. Then the minimal primary decomposition of $I(C_r)$ is
	\[I(C_r)=\bigcap_{m\in M_o(C_r)}(\supp(m))\cap\bigcap_{\{i,i+1\}\in E(C_r)}(s_{i,i+1},t_{i,i+1}).
\]
\end{Proposition}
\begin{proof}``$\supseteq$" By Lemma~\ref{height2and|E|}, $(s_{i,i+1},t_{i,i+1})\supset I(C_r)$ for any $\{i,i+1\}\in E(C_r)$ and they are minimal prime ideals of $I(C_r)$.

Let $m\in M_o(C_r)$ and $\frk{p}=(\supp(m))$. We have to show that $\frk{p}\in\Min(S/I(C_r))$. Firstly, we prove that $\frk{p}\supseteq I(C_r)$. Assume by contradiction $\frk{p}\nsupseteq I(C_r)$, that is there exists a monomial $w\in \mathcal{G}(I(C_r))$ such that $w\notin\frk{p}$. Since, $w\in \mathcal{G}(I(C_r))$, according to Corollary \ref{gencycle}, $w\in M_e(C_r)$. In particular, $\deg_s(w)$ is an even number. On the other hand, $|\supp(m)|=|\supp(w)|=r$ and $\supp(w)\cap\supp(m)=\emptyset$. Thus $\supp(m)\cup\supp(w)=\{s_{i,i+1},t_{i,i+1}\ :\ \{i,,i+1\}\in E(C_r)\}$. Since $m\in M_o(C_r)$, we have $\deg_s(m)$ is an odd number, thus $\deg_t(m)$ is an odd number, which implies $\deg_s(w)$ is odd, a contradiction.

We prove that $\frk{p}\in\Min(S/I(C_r))$. Let $\{\alpha,\alpha+1\}\in E(C_r)$. If $s_{\alpha,\alpha+1}\in\frk p$, let $\frk{p}'=(\supp(m/s_{\alpha,\alpha+1}))$. Assume by contradiction that $\frk{p}'\supseteq I(C_r)$. Since $m\in M_o(C_r)$, $\deg_s(m)$ is odd, which implies $\deg_t(m)$ is an odd number. Let us consider the monomial $M=s_{\alpha,\alpha+1}\prod _{s_{i,i+1}\in\frk{p}'}t_{i,i+1}\prod _{t_{i,i+1}\in\frk{p}'}s_{i,i+1}.$ Since $\deg_t(m)$ is odd, we have $\deg_s(M)$ is an even number. Thus $M\in \mathcal{G}(I(C_r))$, but $M\notin\frk{p}'$, a contradiction.

If $t_{\alpha,\alpha+1}\in\frk p$, let $\frk{p}'=(\supp(m/t_{\alpha,\alpha+1}))$. Assume by contradiction that $\frk{p}'\supseteq I(C_r)$. As before, since $m\in M_o(C_r)$, $\deg_s(m)$ is odd, therefore $\deg_t (m)$ is odd, which implies $\deg_t (m/t_{\alpha,\alpha+1})$ is even. In this case, one may consider the monomial $M=t_{\alpha,\alpha+1}\prod _{s_{i,i+1}\in\frk{p}'}t_{i,i+1}\prod _{t_{i,i+1}\in\frk{p}'}s_{i,i+1}$ which has $\deg_s(M)$ even, thus $M\in \mathcal{G}(I(C_r))$, but $M\notin\frk{p}'$, a contradiction.

``$\subseteq$''  By Lemma~\ref{height2and|E|}, $I(C_r)\subseteq\bigcap _{\{i,i+1\}\in E(C_r)}(s_{i,i+1},t_{i,i+1})$. Let $\frk p\in\Min(S/I(C_r))$,  $\height(\frk p)\geq3$. According to Lemma~\ref{height2and|E|}, we must have $\height(\frk p)=|E(C_r)|$. Let us consider the monomial $m=\prod _{s_{ij}\in\frk p}s_{ij}\prod _{t_{ij}\in\frk p}t_{ij}.$ Assume by contradiction that $\deg_s(m)$ is even, then $\deg_t(m)$ is even. Hence, the monomial $M=\prod _{s_{ij}\in\frk p}t_{ij}\prod _{t_{ij}\in\frk p}s_{ij}$ has $\deg_s(M)$ even, that is $M\in \mathcal{G}(I(C_r))$, but $M\notin\frk p$, a contradiction.
\end{proof}

In order to determine the Castelnuovo--Mumford regularity and the depth of monomial cut ideals of cycles, we need the following preparatory result.

\begin{Lemma}\label{nolinrel} Let $C_r$ be a cycle, $r\geq3$. Then $I(C_r)$ has no linear relations.
\end{Lemma}

\begin{proof} Assume by contradiction that $I(C_r)$ has a linear relation $m_1u_{1}=m_2u_{2}$ with $u_{1},u_{2}\in \mathcal{G}(I(C_r))$ and $\deg(m_1)=\deg(m_2)=1.$ Since $\deg_s(u_{ 1})$ and $\deg_s(u_{ 2})$ are both even, by Proposition~\ref{evensdeg}, we have $m_1,m_2\in\{s_{i,i+1}:1\leq i\leq r\}$ or $m_1,m_2\in\{t_{i,i+1}:1\leq i\leq r\}$. Let us assume that $m_1=s_{i,i+1}$ and $m_2=s_{j,j+1}$, for some $i\neq j$, $1\leq i,j\leq r$, that is $s_{i,i+1}u_{ 1}=s_{j,j+1}u_{ 2}$. This implies $t_{i,i+1}\mid u_{ 1}$, since $s_{i,i+1}\nmid u_{ 1}$ and $\deg(u_{ 1})=r$. Thus $s_{i,i+1}t_{i,i+1}\mid s_{j,j+1}u_{ 2}$ which yields $s_{i,i+1}t_{i,i+1}\mid u_{ 2}$, since $i\neq j$, a contradiction.

If $m_1=t_{i,i+1}$ and $m_2=t_{j,j+1}$, for some $i\neq j$, $1\leq i,j\leq r$, then $t_{i,i+1}u_{ 1}=t_{j,j+1}u_{ 2}$. In particular $s_{i,i+1}\mid u_{ 1}$, since $t_{i,i+1}\nmid u_{ 1}$ and $\deg(u_{ 1})=r$. Thus $s_{i,i+1}t_{i,i+1}\mid s_{j,j+1}u_{ 2}$. But $i\neq j$, hence $s_{i,i+1}t_{i,i+1}\mid u_{ 2}$, a contradiction. 
\end{proof}

As in the case of trees, for monomial cut ideals of cycles, the Castelnuovo--Mumford regularity depends only on the number of edges.

\begin{Proposition}\label{regcycle} Let $C_r$ be a cycle, $r\geq 3$. Then $\reg(S/I(C_r))=|E(C_r)|=r$.
\end{Proposition}

\begin{proof} Let $\Delta$ be the simplicial complex such that $I_{\Delta}=I(C_r)$ and we denote by $I^{\vee}(C_r)$ the Stanley--Reisner ideal of its Alexander dual. By Propositions \ref{primdecodd} and \ref{primdeceven}, $I^{\vee}(C_r)=(M_e(C_r))+(s_{i,i+1}t_{i,i+1}:1\leq i\leq r)$ if $C_r$ is an odd cycle, or $I^{\vee}(C_r)=(M_o(C_r))+(s_{i,i+1}t_{i,i+1}:1\leq i\leq r)$ if $C_r$ is an even cycle. By Lemma~\ref{nolinrel}, we have $\reg(S/I(C_r))\geq r$. According to Terai's Theorem \cite{T}, $\projdim(I^{\vee}(C_r))=\reg(S/I(C_r))$. Hence $\projdim (I^{\vee}(C_r))\geq r$ and $\projdim (S/I^{\vee}(C_r))\geq r+1$. Therefore $\depth (S/I^{\vee}(C_r))\leq r-1$.

 One may note that the $(r-2)$-skeleton of $\Delta^{\vee}$ is $$\Gamma=\Delta^{\vee}_{r-2}=\{F\in\Delta^{\vee}\ :\ |F|=r-1,\ \{s_{i,i+1},t_{i,i+1}\}\nsubseteq F,\ \mbox{for all } 1\leq i\leq r\}.$$ We will show that $\Gamma$ is a shellable simplicial complex, therefore it is Cohen--Macaulay. In particular, $\depth(S/I^{\vee}(C_r))\geq r-1$, which will end the proof. 

First of all, let us fix on $S$ the lexicographical order with $s_{1,2}>\cdots>s_{1,r}>s_{2,3}>\cdots>s_{r,r-1}>t_{1,2}>\cdots>t_{1,r}>t_{2,3}>\cdots>t_{r,r-1}$. For a facet $F_{\alpha}$ of $\Gamma$, we denote by $w_{\alpha}$ the squarefree monomial with the property that $\supp(w_{\alpha})=F_{\alpha}$. We assume that the facets of $\Gamma$ are ordered $F_1,\ldots,F_\ell$ such that, for two facets $F_{\alpha}$ and $F_{\beta}$, $\alpha<\beta$ if $w_{\alpha}>_{lex}w_{\beta}$.

Let $1\leq\beta<\alpha\leq m$. Since $w_{\beta}>_{lex}w_{\alpha}$, there exists $\{i,i+1\}\in E(C_r)$ such that $\{s_{i,i+1},t_{i,i+1}\}\cap(F_{\alpha}\setminus F_{\beta})\neq\emptyset$. If $t_{i,i+1}\in F_{\alpha}\setminus F_{\beta}$, let $F_{\gamma}=F_{\alpha}\cup \{s_{i,i+1}\}\setminus\{t_{i,i+1}\}$ which is a facet of $\Gamma$ since $\{s_{j,j+1},t_{j,j+1}\}\nsubseteq F_{\gamma}$ for any $\{j,j+1\}\in E(C_r)$ and $|F_{\gamma}|=r-1$. We have $w_{\gamma}>w_{\alpha}$, that is $\gamma<\alpha$ and $F_{\alpha}\setminus F_{\gamma}=\{t_{i,i+1}\}$.

Assume that there is no $t_{i,i+1}$ in $F_{\alpha}\setminus F_{\beta}$ for any edge $\{i,i+1\}$. Let $$s_{i,i+1}=\min(F_{\alpha}\setminus F_{\beta})\cap\{s_{l,l+1}:\{l,l+1\}\in E(C_r)\}$$ and $$s_{j,j+1}=\min(F_{\beta}\setminus F_{\alpha})\cap\{s_{l,l+1}:\{l,l+1\}\in E(C_r)\}.$$ Since $w_{\beta}>_{lex}w_{\alpha}$, we must have $j<i$. Let $F_{\gamma}=(F_{\alpha}\setminus\{s_{i,i+1}\})\cup\{s_{j,j+1}\}$. We have $w_{\gamma}>_{lex}w_{\alpha}$ and $F_{\alpha}\setminus F_{\gamma}=\{s_{i,i+1}\}$. One may note that $\{s_{j,j+1},t_{j,j+1}\}\nsubseteq F_{\gamma}$ and for any $l$, $1\leq l\leq r$, $\{s_{l,l+1},t_{l,l+1}\}\nsubseteq F_{\gamma}$, thus $F_{\gamma}$ is a facet of $\Gamma$. Indeed, otherwise we would have $t_{j,j+1}\in F_{\alpha}\setminus F_{\beta}$ which is impossible by our assumption. Therefore, we get that $\Gamma$ is a shellable simplicial complex, hence it is Cohen--Macaulay.

This implies that $\depth(S/I^{\vee}(C_r))=r-1$. Hence $\projdim(S/I^{\vee}(C_r))=r+1$ and, by Terai's theorem, $\reg(I(C_r))=r+1$, that is $\reg(S/I(C_r))=r=|E(C_r)|$.
\end{proof}
Using the above result, one may note that monomial cut ideals of cycles have a pure resolution which is not linear.
\begin{Corollary}\label{betticycle} Let $C_r$ be a cycle, $r\geq3$. Then $\beta_{i,i+j}(S/I(C_r))=0$ for any $j\neq|E(C_r)|$ and $i\geq2$.
\end{Corollary}

\begin{proof} This follows immediately by Proposition \ref{regcycle} and Lemma \ref{nolinrel}.
\end{proof}

We may also compute the depth of monomial cut ideals associated to cycles.

\begin{Proposition}\label{depthcycle} Let $C_r$ be a cycle, $r\geq3$. Then $\depth(S/I(C_r))=|E(C_r)|=r$.
\end{Proposition}
\begin{proof} By Corollary~\ref{betticycle}, $\beta_{i,i+j}(S/I(C_r))=0$ for any $j\neq |E(C_r)|$. Moreover, by Taylor's resolution (cf. \cite[p. 439]{Ei}), for any monomial $m$ in $S$ with $\deg(m)>2|E(C_r)|$, $\beta_{i,\deg(m)}(S/I(C_r))=0$. In particular, $\beta_{i,i+|E(C_r)|}(S/I(C_r))=0$ for any $i> |E(C_r)|$.  This implies $\projdim(S/I(C_r))\leq |E(C_r)|$ and $\depth(S/I(C_r))\geq|E(C_r)|$. By the minimal primary decomposition, we have that $\depth(S/I(C_r))\leq|E(C_r)|$. The statement follows.
\end{proof}
%\begin{Corrolary} Let $G=(V,E(G))$ be a finite connected simple graph. Then $G$ is a cycle if and only if $S/I(G)$ has a pure resolution which is not linear.
%\end{Corrolary}
%\begin{proof} ``$\Rightarrow$" If $G$ is a cycle, the statement follows by Corrolary $$
%\end{proof}

%\begin{Corollary}\label{linres} Let $C_r$ be an cycle, $r\geq3$. Then $\syz_2(S/I(C_r))$ has a linear resolution.
%\end{Corollary}
%\begin{proof} By Corollary~\ref{betticycle}, $S/I(C_r)$ has a pure resolution, therefore $\syz_2(S/I(C_r))$ has a pure resolution. According to Proposition~\ref{depthcycle} 
%\end{proof}
We turn our attention to the minimal primary decomposition and we characterize cycles in terms of the height of minimal prime ideals of monomial cut ideals. 
\begin{Proposition}\label{height} Let $G=(V(G),E(G))$ be a graph. Then $\height(\frk{p})\leq|E(G)|$ for all $\frk{p}\in\Min(S/I(G))$. Moreover, there exists $\frk{p}\in\Min(S/I(G))$ of height $|E(G)|$ if and only if $G$ is a cycle or $G$ is the path graph $ 1-2-3$. 
\end{Proposition}
\begin{proof} First of all, one may note that $(s_{ij},t_{ij})\supset I(G)$ and they are minimal prime ideals. If we assume by contradiction that there exists $\frk{p}\in\Min(S/I(G))$ such that $\height(\frk p)>|E(G)|$, then there exists $\{i,j\}\in E(G)$ such that $(s_{ij},t_{ij})\subset \frk p$, contradiction with $\frk p\in\Min(S/I(G))$.

``$\Leftarrow$" This obviously follows by the minimal primary decomposition of monomial cut ideals of cycles and trees.

``$\Rightarrow$" Let $\frk{p}\in\Min(S/I(G))$ such that $\height(\frk{p})=|E(G)|$. 
If $G$ is a tree, then, by Proposition~\ref{primdectree}, $\Min(S/I(G))=\{(s_{ij},t_{ij}):\{i,j\}\in E(G)\}$, which implies $\height(\frk{p})=2$, for all $\frk{p}\in\Min(S/I(G))$. Therefore $G$ is the graph path graph $ 1-2-3$. Next, we will assume that $G$ is not a tree. 

Assume by contradiction that $G$ is not a cycle. Therefore, $G$ contains a proper subgraph which is a cycle. Let $C_r\subsetneq G$ be a cycle and $\{\alpha,\beta\}\in E(C_r)$. Since $\height(\frk{p})=|E(G)|$ and $\frk{p}\in\Min(S/I(G))$, we see that, for every $\{i,j\}\in E(G)$, $s_{ij}\in\frk{p}$ or $t_{ij}\in\frk{p}$. 

If $s_{\alpha\beta}\in\frk{p}$, let $\frk{p}=\frk{p}'+(s_{\alpha\beta})$ and $(s_{\alpha\beta})\nsubseteq\frk{p}'$. Let $\frk{p}_r=\frk{p}\cap\kk[s_{ij},t_{ij}:\{i,j\}\in E(C_r)]$ and $u_{\frk{p}_r}$ be the squarefree monomial such that $(\supp(u_{\frk{p}_r}))=\frk{p}_r$. If we prove that $\frk{p}'\supset I(G)$, the statement follows. 

Assume by contradiction that $\frk{p}'\nsupseteq I(G)$, that is there exists $m\in \mathcal{G}(I(G))$ such that $m\notin\frk{p}'$. One may note that we cannot have $m=\prod_{s_{ij}\in\frk{p}}t_{ij}\prod_{t_{ij}\in\frk{p}}s_{ij},$ since $m\in I(G)\subseteq \frk{p}$, thus we must have $m=s_{\alpha\beta}\prod_{s_{ij}\in\frk{p}'}t_{ij}\prod_{t_{ij}\in\frk{p}'}s_{ij}\in I(G).$ Therefore, we may consider $m_r=s_{\alpha\beta}\prod_{s_{ij}\in \frk{p}_r,\ \{i,j\}\neq\{\alpha,\beta\}}t_{ij}\prod_{t_{ij}\in \frk{p}_r}s_{ij}\in I(C_r)$ and we have $\deg_s(m_r)$ is an even number.

If $C_r$ is an even cycle, then $\deg_t(m_r)$ is an even number and we obtain that $\deg_s(u_{ \frk{p}_r})=\deg_t(m_r)+1$ is an odd number. Hence $ \frk{p}_r\in\Min(S/I(C_r))$ and $\frk{p}\supsetneq \frk{p}_r\supset I(C_r)\supset I(G)$, a contradiction.

If $C_r$ is an odd cycle, then $\deg_t(m_r)$ is an odd number and $\deg_s(u_{ \frk{p}_r})=\deg_t(m_r)+1$ is an even number. Hence $ \frk{p}_r\in\Min(S/I(C_r))$ and $\frk{p}\supsetneq \frk{p}_r\supset I(C_r)\supset I(G)$, a contradiction.

The same reasoning applies to the case when $t_{\alpha\beta}\in\frk{p}$. We write $\frk{p}=\frk{p}'+(t_{\alpha\beta})$ and $(t_{\alpha\beta})\nsubseteq\frk{p}'$, and we consider the ideal $ \frk{p}_r=\frk{p}\cap\kk[s_{ij},t_{ij}:\{i,j\}\in E(C_r)]$. Let $u_{ \frk{p}_r}$ be the squarefree monomial such that $(\supp(u_{ \frk{p}_r}))= \frk{p}_r$. We will show that $\frk{p}'\supset I(G)$ which will end the proof. 

Assume by contradiction that $\frk{p}'\nsupseteq I(G)$, thus there exists $m\in \mathcal{G}(I(G))$ such that $m\notin\frk{p}'$. As before, since $m\in I(G)\subseteq \frk{p},$ we must have $m=t_{\alpha\beta}\prod_{s_{ij}\in\frk{p}'}t_{ij}\prod_{t_{ij}\in\frk{p}'}s_{ij}\in I(G).$ Therefore, we may consider $m_r=t_{\alpha\beta}\prod_{s_{ij}\in \frk{p}_r}t_{ij}\prod_{t_{ij}\in \frk{p}_r,\ \{i,j\}\neq\{\alpha,\beta\}}s_{ij}\in I(C_r)$ and we have $\deg_s(m_r)$ is an even number.

If $C_r$ is an even cycle, then $\deg_t(m_r)$ is an even number and we obtain that $\deg_s(u_{ \frk{p}_r})=\deg_t(m_r)-1$ is an odd number, which implies that $ \frk{p}_r\in\Min(S/I(C_r))$ and $\frk{p}\supsetneq \frk{p}_r\supset I(C_r)\supset I(G)$, a contradiction.

If $C_r$ is an odd cycle, then $\deg_t(m_r)$ is an odd number and $\deg_s(u_{ \frk{p}_r})=\deg_t(m_r)-1$ is an even number. Thus $ \frk{p}_r\in\Min(S/I(C_r))$ and $\frk{p}\supsetneq \frk{p}_r\supset I(C_r)\supset I(G)$, a contradiction.
\end{proof}
\section{Primary decomposition of monomial cut ideals}

We aim at describing the minimal primary decomposition of monomial cut ideals associated to an arbitrary graph. First of all, we describe the relation between the monomial cut ideal of a graph and the monomial cut ideals of its proper subgraphs.

\begin{Proposition}\label{subgraph1} Let $G=(V(G),E(G))$ be a graph which is not a cycle. Then $$I(G)=\bigcap_{\twoline{H\mbox{\rm\scriptsize{ a proper}}}{\mbox{\rm\scriptsize{ subgraph of }} G}}I(H).$$
\end{Proposition}

\begin{proof} ``$\supseteq$" Let $H$ be a proper subgraph of $G$ and $\frk{p}\in\Min(S/I(H))$. Hence $\frk{p}\supseteq I(H)\supseteq I(G)\supseteq I(H)\cdot\left(\prod _{\{i,j\}\in E(G)\setminus E(H)}s_{ij}t_{ij}\right).$ Since $\prod _{\{i,j\}\in E(G)\setminus E(H)}s_{ij}t_{ij}\notin\frk{p}$, we have $\frk{p}\in\Min(S/I(G))$.

``$\subseteq$" If $G$ is the path graph $1-2-3$, then the statement follows by Proposition~\ref{primdectree}. Therefore, we may assume that $G$ has at least three edges and $G$ is not a cycle.

Let $\frk{p}\in\Min(S/I(G))$. By Proposition~\ref{height}, $\height(\frk{p})<|E(G)|$. Therefore, there exists an edge $\{\alpha,\beta\}\in E(G)$ such that $s_{\alpha\beta}t_{\alpha\beta}\notin\frk{p}$. Let $H$ be the subgraph of $G$ obtained by deleting the edge $\{\alpha,\beta\}$. By Proposition~\ref{edge}, $I(H)=I(G):(s_{\alpha\beta}t_{\alpha\beta})$, therefore, $\frk{p}\supseteq I(G):(s_{\alpha\beta}t_{\alpha\beta})=I(H)$. Using Remark~\ref{subgraph}, we get $$\frk{p}\supseteq I(H)\supseteq\bigcap_{\twoline{\tilde{H}\mbox{\rm\scriptsize{ a proper}}}{\mbox{\rm\scriptsize{ subgraph of }} G}}I(\tilde{H})\supseteq I(G),$$ and the statement follows.
\end{proof}

\begin{Theorem}\label{primdec} Let $G=(V(G),E(G))$ be a graph. Then the minimal primary decomposition of $I(G)$ is
\[I(G)=\bigcap_{\twoline{C\mbox{\rm\scriptsize{ a cycle of }} G}{\twoline{{\bf{p}}\in\Min(I(C))}{\height({\bf{p}})>2}}}{\bf{p}}\cap\bigcap_{\{i,j\}\in E(G)}(s_{ij},t_{ij}).
\]
\end{Theorem}

\begin{proof} We use induction on the number of edges. If $G$ has only one edge $\{1,2\}$, then $I(G)=(s_{12},t_{12})$. Assume that the statement is true for any graph with $r\geq1$ edges and let $G$ be a graph with $r+1$ edges. The statement is true for trees and cycles, by Propositions~\ref{primdectree}, \ref{primdecodd}, and \ref{primdeceven}. Assume that $G$ is neither a tree nor a cycle, thus $G$ contains at least a cycle as a proper subgraph. Then, by Proposition~\ref{subgraph1} and by the induction hypothesis, the statement follows.
\end{proof}
\begin{Corollary}\label{Krulldim}If $G=(V(G),E(G))$ is a graph, then $\dim(S/I(G))=2|E(G)|-2$.
\end{Corollary}
We may determine the multiplicity of monomial cut ideals associated to an arbitrary graph.

\begin{Corollary}\label{multiplicity} Let $G=(V(G),E(G))$ be a graph. Then $e(S/I(G))=|E(G)|.$
\end{Corollary}
\begin{proof} Let $\Delta$ be the simplicial complex such that $I_{\Delta}=I(G)$. By Corollary~\ref{Krulldim}, $\dim(\Delta)=2|E(G)|-3$. Since we have precisely $|E(G)|$ faces of dimension $2|E(G)|-3$, we get that $e(S/I(G))=f_{2|E(G)|-3}=|E(G)|$. 
\end{proof}

We are now interested in relating algebraic properties of monomial cut ideals such as being unmixed, or Cohen--Macaulay with combinatorial properties of the graph. First of all, we characterize all monomial cut ideals which are unmixed.
 
\begin{Theorem} Let $G=(V(G),E(G))$ be a graph and $I(G)$ the corresponding monomial cut ideal. The following statements are equivalent:
\begin{itemize}
\item[a)] $I(G)$ is unmixed.
\item[b)] $I(G)$ has a linear resolution.
\item[c)] $I(G)$ has only linear relations.
\item[d)] $G$ is a tree.
\end{itemize}
\end{Theorem}
\begin{proof} One may note that ``(a)$\Leftrightarrow$(d)" follows by Proposition~\ref{primdectree} and Theorem~\ref{primdec}.

``(a)$\Rightarrow$(b)" Let $\Delta$ be the simplicial complex such that $I_{\Delta}=I(G)$. By Theorem~\ref{primdec}, $I_{\Delta^{\vee}}=(s_{ij}t_{ij}:\{i,j\}\in E(G))$ which is a complete intersection ideal, thus it is Cohen--Macaulay. The statement follows by the Eagon--Reiner Theorem \cite{EaRe}.

The implication ``(b)$\Rightarrow$(c)" obviously holds.

``(c)$\Rightarrow$(d)" Assume by contradiction that $G$ contains a cycle $C$. Then, any relation of $I(C)$ gives us a relation of $I(G)$. According to Lemma~\ref{nolinrel}, $I(C)$ has no linear relations, a contradiction. 
\end{proof}
Concerning the property of being Cohen--Macaulay, we have the following characterization:
\begin{Theorem} Let $G=(V(G),E(G))$ be a graph and $I(G)$ the corresponding monomial cut ideal. The following statements are equivalent: 
\begin{itemize}
\item[a)] $I(G)$ is a complete intersection ideal.
\item[b)] $I(G)$ is Gorenstein.
\item[c)] $I(G)$ is Cohen--Macaulay.
\item[d)] $I(G)$ is sequentially Cohen--Macaulay.
\item[e)] $G$ is the graph with only one edge.
\end{itemize}
\end{Theorem}
\begin{proof} The implications ``(a)$\Rightarrow$(b)$\Rightarrow$(c)$\Rightarrow$(d)" are known, and one may easy see that ``(e)$\Rightarrow$(a)" is also true. We have to prove that ``(d)$\Rightarrow$(e)". Let $\Delta$ be a simplicial complex such that $I_{\Delta}=I(G)$. By \cite[Theorem 2.1(a)]{HH}, $I(G)$ is sequentially Cohen--Macaulay if and only if $I_{\Delta^{\vee}}$ is componentwise linear. In particular, $(I_{\Delta^{\vee}})_{\langle2\rangle}$ has a linear resolution, where $(I_{\Delta^{\vee}})_{\langle2\rangle}$ is the monomial ideal generated by all the monomials of degree 2 from $I_{\Delta^{\vee}}$. By Theorem~\ref{primdec}, $(I_{\Delta^{\vee}})_{\langle2\rangle}=(s_{ij}t_{ij}:\{i,j\}\in E(G))$, therefore $(I_{\Delta^{\vee}})_{\langle2\rangle}$ is a complete intersection ideal generated in degree $2$, which has a linear resolution if and only if it is a principal ideal. This implies that $|E(G)|=1$.
\end{proof}
	
\begin{center}

\textbf{Acknowledgments}

\end{center}

Part of this paper was prepared during the visit to the Duisburg--Essen University, Germany. The author would like to thank Professor J\"urgen Herzog for the kind hospitality, for useful discussions and valuable suggestions.

\end{document}